\documentclass[12pt]{iopart}

\renewcommand{\vec}{\bi}
\newcommand{\E}{\mbox{I\negthinspace E}}

\usepackage{iopams}
\usepackage{graphicx}
\usepackage{dsfont}

\begin{document}

\title{Optimal control theory with arbitrary superpositions of waveforms}

\author{Selina Meister$^1$, J\"urgen T. Stockburger$^1$, Rebecca Schmidt$^2$ and Joachim Ankerhold$^1$}
\address{$^1$\ Institut f\"ur komplexe Quantensysteme, Universit\"at Ulm,
Albert-Einstein-Allee 11, D-89069 Ulm, Germany}
\address{$^2$\ School of Mathematical Sciences, The University of Nottingham, Nottingham, NG7 2RD, United Kingdom}


\date{\today}

\begin{abstract}
Standard optimal control methods perform optimization in the time
domain.  However, many experimental settings demand the expression
of the control signal as a superposition of given waveforms, a case
that cannot easily be accommodated using time-local constraints.
Previous approaches \cite{Skinner10,Moore12} have circumvented this 
difficulty by performing optimization in a parameter space, using the
chain rule to make a connection to the time domain.
In this paper, we present an extension to Optimal Control Theory which
allows gradient-based optimization for superpositions of arbitrary
waveforms directly in a time-domain subspace. Its key is the use of the
Moore-Penrose pseudoinverse as an efficient means of transforming between
a time-local and waveform-based descriptions. To illustrate this
optimization technique, we study the parametrically driven harmonic
oscillator as model system and reduce its energy, considering both
Hamiltonian dynamics and stochastic dynamics under the influence of a
thermal reservoir. We demonstrate the viability and efficiency of the
method for these test cases and find significant advantages in the
case of waveforms which do not form an orthogonal basis.
\end{abstract}

%
%
\submitto{\jpa}
%
\maketitle
%
%
\section{Introduction}
Optimal control theory aims at driving a dynamical system towards a
final state that minimizes a figure of merit and at finding the
required time-dependent controls. In classical physical systems,
optimal control schemes have been used successfully for decades
\cite{Bryson75,Davis2002}, and for some time various optimal control
algorithms have been applied to a wide range of quantum systems, see
e.g.\ \cite{Brif10
  ,Glaser1998
  ,Sklarz2002,palao02,Calarco2000,Schmidt13}. 
There is a renewed interest in applying and improving established
algorithms like the Krotov algorithm \cite{Krotov96} or gradient
methods \cite{Khaneja05}, and new control techniques are being
developed for control problems of increasing complexity
\cite{Haeberle14,Goerz14,
  deFouquieres11, Caneva11}.

Speed-up of optimization is often achieved by restriction of the
control pulses to certain pulse families, such as, e.g., Gaussian
pulse cascades \cite{Emsley90} or Fourier expansions \cite{Bartels13,
  RomeroIsart07, Doria11,Verdeny14
}. Also, more elaborate ways of truncating the search space have been
formulated, like e.g. the CRAB algorithm \cite{Caneva11}.
These methods are especially successful, if the underlying dynamics is
well known and understood, which allow for a sophisticated choice of
the basis functions. Additionally, the comparatively easy shape of the
pulse (limited, e.g., to a few frequency components only), allows for
a straight forward interpretation.

In order to gain experimentally realizable control pulses, additional
constraints may have to be taken into account, such as restricting the
total energy or limiting the control functions
\cite{Khaneja05,Haeberle14, Graichen09,
  Graichen10,Schaefer12,Sugny14,
  Palao13
}. Not all control algorithms follow this requirement equally well;
for the effects such constraints have on the convergence behaviour of
control algorithms, see \cite{Moore12}. Methods to design pulse shapes as analytic functions of a small set of parameters have been
introduced \cite{Skinner10}, which allow restrictions on the shape of
the control pulse. The most elementary extension of gradient-based
control theory to this scenario consists in simply applying the chain
rule of calculus.

In this paper we present a control algorithm which takes similar
experimental constraints into account, while working in a proper
subspace of actual functions of time rather than a parameter space.
Our control subspace is defined from linear combinations of arbitrary
waveforms. The properties of orthogonality, normalization, or even
linear independence are not required in our case, valid solutions in
the time domain are obtained in either case.

The method is tested by examining a generic model system, for which
methods of comparison are available. However, the method is not
restricted to such simple systems, but, as in section \ref{sec:method} described,
broadly applicable.

The paper is organized as follows: In section~\ref{sec:method} the
general method is described which is then applied to generic test
cases in section~\ref{sec:systemandresults}, including a discussion of
performance characteristics. A short summary together with an outline
of potential extensions and applications  follow in section \ref{sec:conclusion}.

\section{Method: experimentally realizable control
  functions\label{sec:method}}

The aim is to control a dynamical system with an equation of motion
$\dot{\vec{z}}=f(\vec{z}(t),u(t))$, where $z(t)$ is the dynamical
state vector, and $u(t)$ is an external time-dependent control. An
optimal control function $u(t)$ is sought which changes the dynamics
of the system towards a desired property of a final state at given
time $\tau$. The objective can be quantified by a cost functional $
\Phi[\vec{z}(\tau)]$. For the sake of clarity, only a scalar control
signal $u(t)$ is considered; but this is easily generalized to the
case of multiple control variables.  To find the optimal control
function $u(t)$, which also considers the constraint that the
equations of motion needs to be satisfied, an extremum of the
augmented cost functional
\begin{equation}
\label{eq:J}
 J[u(t)]=\Phi[\vec{z}(\tau)]+\int\limits_{t_i}^\tau \vec{\lambda}^t(t)\{\dot{\vec{z}}(t)-f(\vec{z}(t),u(t))\} dt ,
\end{equation}
needs to be determined. $\vec{\lambda}(t)$ denotes a vector of
Lagrange multipliers.

Variational calculus on (\ref{eq:J}) yields equations of
motion also for the Lagrange multipliers, which are therefore referred
to as co-states. Variation of the system state at the final time
yields a final-time boundary condition for the co-state, which depends
on the final-time value of the state variable
$\vec{z}(\tau)$. Furthermore, variation with respect to $u(t)$ leads to
the gradient
\begin{equation}\label{eq:grad}
 \quad \frac{\partial J}{\partial u(t)}=-\vec{\lambda}^t(t)\frac{\partial f}{\partial u}.
\end{equation}
The optimal control function $u(t)$ can then be found by searching for zeros of this gradient.
For solving optimal control problems gradient methods are well established \cite{Bryson75,Khaneja05}. As computing the gradient is relatively cheap (only twice the effort of a simple evaluation of the cost functional), it is almost always favourable to adopt an optimization technique which makes use of the information supplied by the gradient to determine the direction and width of the search step for the next optimization iteration. The method presented here is therefore based on a gradient method.
\begin{figure}
\begin{center}
\includegraphics[width=10cm]{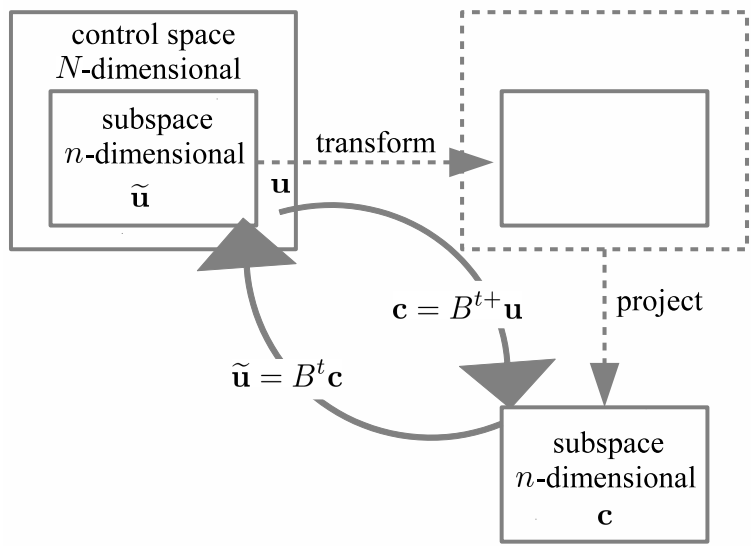} 
\caption{Schematic representation of the function
  spaces in which the optimization problem can be formulated. In full
  control space (big box in the upper left) the time discretized
  control pulse is a $N$-dimensional vector $\vec{u}$. An
  $n$-dimensional subspace (inset, top left), spanned by vectors
  combined in a matrix $B$, the control pulse may be described by a
  $n$-dimensional coefficient vector $\vec{c}$ (bottom). After
  transforming to the subspace and back, $\vec{u}$ is not recovered,
  but $\vec{\widetilde{u}}$, a time dependent vector in the
  $n$-dimensional subspace. Note that a direct transformation (curved
  down arrow) based on the Moore-Penrose pseudoinverse $B^+$ performs
  both a projection and a change of representation in one step.}
\label{fig:spaces}
\end{center}
\end{figure}

 The numerical representation of the control function $u(t)$ is, due to the discrete time steps $\delta t$, a $N$-dimensional vector:
\begin{equation}
u(t)\rightarrow (u(t_1),\ldots ,u(t_k), \ldots, u(t_N)), \quad
  t_k=k\cdot\delta t, \quad k=1,\ldots ,N, 
\end{equation}
with $\tau\equiv t_N$, $t_i\equiv t_1$ and in shorter notation
$u(t)\rightarrow (u_1,\ldots ,u_k, \ldots, u_N) \equiv
\vec{u}^t$. With the dynamics as an implied constraint, the cost to be
optimized becomes a functional $\Phi[\vec{u}]$.  The underlying
$N$-dimensional function space will be called full function space or
control space.

However, depending on the problem under study, additional constraints
for the control function $u(t)$, the dynamical variable $\vec{z}(t)$
or a function $\mathcal{E}(\vec{z},u,t)$ may be required. Hard
boundaries for the variables just mentioned can be expressed as
inequality conditions $\mathcal{E}(\vec{z},u,t)>0$. There are
different approaches to consider these restrictions: Either additional
Lagrange terms in the augmented cost functional (\ref{eq:J}) are
included in the formalism or a suitable change of variables is
performed. For details see  \cite{Haeberle14, Graichen09,
  Graichen10}. Another way is to take into account additional cost
functionals of the form $J_{add}[u(t)]=\alpha\int\limits_{t_i}^\tau
dt\Theta[\mathcal{E}(\vec{z},u,t)]$, where $\Theta$ denotes the
Heaviside-function and $\alpha$ is a parameter to give priority to this
constraint. A large enough value of $\alpha$ enforces the restriction,
but there is no guarantee that it will not be violated for brief
periods of time, when using additional cost functionals.

On the other hand, an important, experimentally relevant class of
constraints cannot be expressed through inequality constraints or an
additional time-local term in the cost functional. E.g, pulse shapes
available for optical control signals are often characterized by a set
of parameters determining a pulse sequence through central frequencies
and shape parameters. Such constraints could be expressed through
equations involving functionals of the control signal, however, we are
going to take the simpler, more transparent approach of working in the
lower-dimensional (truncated) function space implied by the
parametrization. Note that this is subtly different from working
directly in the parameter space~\cite{Skinner10}: Even if both spaces
share the same affine structure, their natural scalar products, which
are used in numerical optimization, are rarely identical.

\begin{figure}
\begin{center}
\includegraphics[width=10cm]{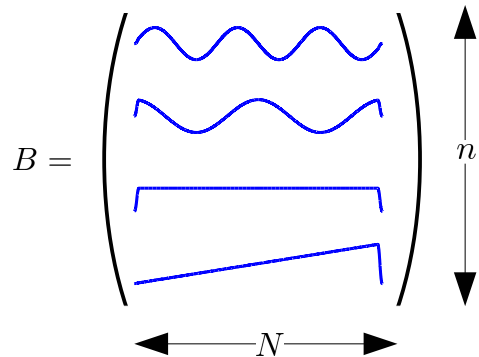} 
\caption{(color online) Schematic representation of the Matrix
$B$ (\ref{eq:B}) with rows determined by the waveforms $b_l(t)$ or
their vector equivalents $\vec{b}_l^t$. The transpose $B^t$ transforms a
coefficient vector $\vec{c}$ into a control function $\vec{u}$, see
equation (\ref{eq:ctou}). $B$ is not square, typically $n\ll N$.
The vectors $\vec{b}_l$ representing the waveforms need not
to be orthogonal.}\label{fig:basis}
\end{center}
\end{figure}

In this paper we present a method to take into account the type of
constraint implied by a parametrization of the control $u(t)$ as the
superposition of $n$ experimentally realizable waveforms $b_l(t)$ (cf. figure \ref{fig:spaces}).
The functions $b_l(t)$ can be interpreted as vectors, defining a linear space through their span. This space is obviously a subspace of the space of control functions $u(t)$ admissible in the absence of constraints.
In the numerical representation, we get $N$-dimensional vectors
$\vec{b}_l$. These vectors form an $n\times N$ matrix
\begin{equation}
  \label{eq:B}
B=
\left(\begin{array}{c}
  \vec{b}_1^t \\ \vdots \\ \vec{b}_n^t
\end{array}\right),
\end{equation}
which can be used to tranform the coefficient vector $\{ c_1,\ldots
c_l,\ldots ,c_n \}\equiv \vec{c}^t$ of a superposition of waveforms into
the time domain, i.e., the superposition
\begin{equation}
 \widetilde{u}(t)=\sum\limits_{l=1}^n c_lb_l(t)
\end{equation}
can be written in the compact form
\begin{equation}
\label{eq:ctou}
 \vec{\widetilde{u}} =  B^t\vec{c}.
\end{equation}
The matrix  $B$ is schematically shown in figure \ref{fig:basis}.

The optimization now consists in the task of finding a control
$\vec{\widetilde{u}}$ within the subspace for which
$\Phi[\vec{\widetilde{u}}]$ is at an extremum.  We shall see below
that this is feasible without placing any restrictions on the vectors
$\{\vec{b}_l\}$. Note that the $\{\vec{b}_l\}$ need not be linearly
independent, in particular they need not be normalized or orthogonal,
but we may choose them normalized for convenience.

Lack of any restrictions on the vectors $\{\vec{b}_l\}$ makes our
approach highly flexible. Functions $\{b_l(t)\}$ with arbitrary shapes
may be chosen - symmetric or asymmetric pulses, ramps and plateaus,
chirped pulses, virtually anything a given application may
require. Since orthogonality is not a criterion, pulse shapes can also
easily be modified for smooth rise or fall to/from an initial or final
value of zero.

Now we want to minimize $\Phi[\vec{\widetilde{u}}]$ with the aid of
the gradient (\ref{eq:grad}). The gradient can be calculated by first
solving the equation of motion for the system degrees of freedom
$\dot{\vec{z}}=f(\vec{z},\widetilde{u})$. Then the final time value
$\vec{z}(\tau)$ fixes the end-time value of the Lagrange multiplier,
which can be propagated backwards in time to get the solution
$\vec{\lambda}(t)$. This is the standard computational approach to
obtain the gradient (\ref{eq:grad}). In the present context, the
equations of motion and the gradient depend on the truncated control
$\vec{\widetilde{u}}$, as the optimization is performed in the
subspace. However, the gradient vector obtained from (\ref{eq:grad})
typically lies outside the subspace of interest, it needs to be
projected back each time it is computed. It is therefore
essential to find an efficient way to perform this projection as well
as the transformations between the time domain and the coefficients of
the waveform decomposition.

Since the waveforms cannot be assumed to be orthogonal functions, the
transformation cannot be performed using a simple scalar product
\begin{equation}
c_l=\vec{u}\cdot \vec{b}_l \quad \not\Leftrightarrow \quad \vec{\widetilde{u}}=\sum\limits_{l=1}^nc_l\vec{b}_l.
\end{equation}
This also means that $\sum\limits_l \vec{b}_l^t\cdot\vec{b}_l$ cannot
be used as a projector on the constrained subspace, which would be the
case for an orthonormal basis. Moreover, if the $\{\vec{b}_l\}$ are
not linearly independent, the coefficients $c_l$ are not
unique. However, we shall demonstrate below that even in this case the
constrained optimization problem is well-defined and converges to
unambiguous solutions $\widetilde{u}(t)$.

Even for the case of linearly independent $\{\vec{b}_l\}$, completing
the matrix $B$ into a quadratic matrix with full rank, computing its
inverse and performing a projection onto the subspace (dotted indirect
path in \ref{fig:spaces}) would generally be much more tedious than the method
we introduce here.

To be able to perform the transformation with ease, we suggest the use of the Moore-Penrose-Pseudoinverse $B^+$ \cite{Penrose55} (PINV). It has the following defining properties:
\begin{equation}
\eqalign{ BB^+ B =&B  \qquad \quad \; B^+ BB^+ =B^+ \cr
  (BB^+)^* =& BB^+ \qquad (B^+ B)^* = B^+ B. }
\end{equation}
$B^+$ exists for \emph{any} matrix $B$. Both $B^+B$ and $BB^+$ are projection matrices. This definition contains the ordinary inverse as a special case, where $B^+B=BB^+=\mathds{1}$. The equations
\begin{equation}
\vec{c}=B^{t+}\vec{u} \quad \mbox{and} \quad
\vec{\widetilde{u}}=B^t\vec{c} = B^t{B^{t+}} \vec{u} = B^+B \vec{u}
\end{equation}
provide all the required transformations and projections (see figure
\ref{fig:spaces}).

When transforming back to the time domain in the case $n<N$, it is not
necessarily the original vector $\vec{u}$ that is recovered, but
$\vec{\widetilde{u}}$, which is the result of the projection
$B^+B$. This is exactly what is required for a consistent
iterative optimization in the subspace.
Suitable algorithms for the precise and efficient
computation of the pseudoinverse exist and are implemented in many
numerical libraries and software packages.

We are now in a position to formulate an iteration scheme which
computes a solution of the \emph{constrained} control problem while
making use of the gradient of the objective functional in the
\emph{full} space, which can be obtained from the dynamics of states
and co-states. Now the following step-by-step instruction for the PINV method can be given:
\begin{enumerate}
\item Define the cost functional and choose the generating system $\{\vec{b}_l\}$ for the subspace. Choose an initial guess $\vec{u_i}$ for the optimal control pulse.
\item Compute the values of the cost functional and the gradient and project the gradient on the subspace:
\begin{equation}
J[\vec{\widetilde{u}}_i]=J[B^+B\vec{u}_i]
\end{equation}
\begin{equation}
\nabla_{\widetilde{u}}J[\vec{\widetilde{u}}_i]=B^+B\,
\nabla_uJ[\vec{u}]|_{\vec{u}=\vec{\widetilde{u}}_i}.
\end{equation}
$\nabla_u$ denotes the gradient with respect to $\vec{u}$.
\item Perform the optimization search step, e.g., a line search, in
  the $n$-dimensional subspace. This can
  be done by a ``black box'' routine implementing standard
  optimization techniques. In most cases, one thus finds a control signal
  $\vec{\widetilde{u}}_{i+1}$ with
  $J[\vec{\widetilde{u}}_{i+1}]<J[\vec{\widetilde{u}}_{i}]$.
\item Set $\vec{u}_i=\vec{\widetilde{u}}_{i+1}$ and start again at
  step (ii) until
  $J[\vec{\widetilde{u}}_{i}]-J[\vec{\widetilde{u}}_{i+1}]<\epsilon$, where
  the parameter $\epsilon$ then defines the convergence criterion.
\end{enumerate}

Using these steps, both the gradient search and conjugate gradient
algorithms as well as quasi-Newton optimization algorithms
\cite{Boyd04, Fletcher87, Nocedal06} can be adapted to perform
optimal control in a parameter subspace. The latter algorithms use
values of the gradient from several iteration cycles to build an
approximate Hessian matrix to accelerate convergence.

\section{Test cases and results \label{sec:systemandresults}}
As a proof of concept, we apply our method to two simple model
systems, for which comparison data can easily be obtained. In order to
show the flexibility of the tested PINV method and the compatibility
with different optimization algorithms, the results presented here
are computed using different gradient-based optimization
methods. In section \ref{sec:deterministic} a L-BFGS
  quasi-Newton method \cite{Nocedal06} and in section \ref{sec:open} a
  steepest descent algorithm will be used.

\subsection{Control of deterministic dynamics \label{sec:deterministic}}

\begin{figure}
\begin{center}
  \includegraphics[width=10cm]{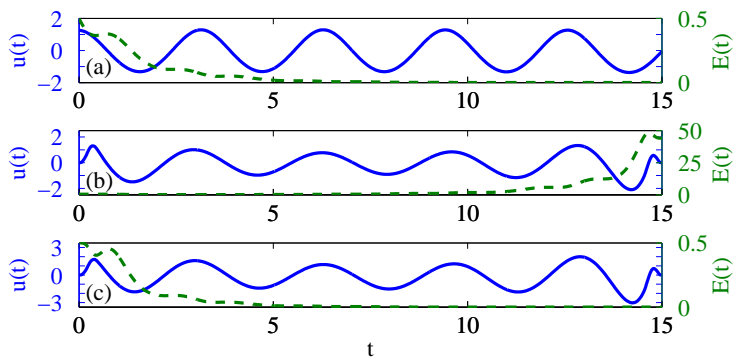}\caption{(color online) The control signals
    (blue solid lines, dimensionless units on left axes) and the
    energies of the controlled system (green dashed lines, right
    axes). Data in (a) correspond to the optimization in full function
    space with final-time energy
    $E(\tau)=2.850\cdot 10^{-5}$. In (b) this control is truncated to
    a $n=12$-dimensional subspace defined in (\ref{eq:basisp}). This
    post-truncated control leads to $E(\tau)=44.096$. In (c) the
    result achieved with the PINV method in the same subspace as in
    (b) is shown. Using this control function one achieves
    $E(\tau)=4.432\cdot 10^{-6}$}.\label{fig:T0}
\end{center}
\end{figure}

The system of interest for this section is a parametrically driven
harmonic oscillator. Its dynamics is described by the set of
differential equations
\begin{equation} \label{eq:HOT0}
\eqalign{\dot{q}&=\frac{p}{m} \cr
\dot{p}&=-m\omega_0^2q(t)-u(t)q(t)}
\end{equation}
with a parametric control function $u(t)$. For convenience, we set
$\omega_0=1$ and $m=1$ in the following (natural units) and choose an
initial state $(q_0,p_0)=(1/\sqrt{2},1/\sqrt{2})$ with initial energy
$E(0)=0.5$. The control objective is the minimization of energy at the
fixed final time $\tau=15$, i.e.
\begin{equation}
\Phi[q(\tau),p(\tau)]=\frac{q^2(\tau)}{2}+\frac{p^2(\tau)}{2}.
\end{equation}
As generating set for the subspace in which the optimization will be performed we choose the following (non-orthogonal) functions
\begin{equation}
\left.
\begin{array}{c}
b_l(t) =\frac{1}{\mathcal{N}}e(t)\sin\left(\frac{2\pi}{\tau}l\cdot t\right) \\\\
b_{l+\frac{n}{3}}(t) =\frac{1}{\mathcal{N}}e(t)\cos\left(\frac{2\pi}{\tau}l\cdot t\right)\\\\
b_{l+\frac{2n}{3}}(t) =\frac{1}{\mathcal{N}}e(t)t^{l-1}
\end{array}  \right\} \quad l=1,\cdots, \frac{n}{3}
\label{eq:basisp}
\end{equation}
on $[0,\tau]$ with a scaling constant $\mathcal{N}$ and an envelope
function $e(t)$ enforcing initial and final values $b_k(0) = b_k(\tau)
= 0$, modeled on the finite time needed to switch a real-world pulse.
The function $e(t)$ rises from zero to unity within a specified time
interval $[0,t_0]$, followed by a plateau and a symmetric decline in
the interval $[\tau-t_0,\tau]$. For the initial rise, $t\in [0,t_0]$,
we choose
\begin{equation}\label{eq:envelope}
e(t)=\frac{1}{2}\left\{ 1+\tanh\left[\frac{\eta \cdot y(t)}{1- y(t)^4}  \right]  \right\},
\end{equation}
with $y(t) =\frac{2t}{t_0}-1$ and the parameters $t_0 = 0.5$ and
$\eta=2$. Note that $e(t)$ is smooth at $t=t_0$ and $t=\tau-t_0$, all
its derivatives are zero at these points. The
matrix $BB^t$, which is the unit matrix in the case of an orthonormal
basis, has widely differing eigenvalues in our case; its condition
number is approximately $1.3\cdot 10^8$.

In figure \ref{fig:T0} the control functions (blue solid lines) and
the corresponding energies (green dashed lines) are shown for three
different approaches to optimal control. The top panel displays
results for optimization in the full space.
Using the L-BFGS quasi-Newton algorithm \cite{Nocedal06}, a solution
with a very small final energy $E(\tau)=2.850\cdot 10^{-5}$ is
computed. (Since this is a fully controllable system, the exact solution is zero).
As might be expected for parametric control, the dominant
frequency of this control is twice the system frequency $2\omega_0$,
which can be qualitatively understood as follows: When the particle is
in the minimum of the potential, the control makes the potential well
narrower, decreasing its amplitude after traversing the minimum. When
the particle is near the turning points, the potential is opened to
extract energy.

In panel (b) of figure \ref{fig:T0}, the control pulse optimized in full function
space [panel (a)] is truncated to the subspace after the
iteration has been performed, instead of using the PINV scheme
described in section \ref{sec:method}. The basis functions are given by (\ref{eq:basisp}) with $n=12$. Although the truncated control
signal still looks fairly similar to the original one, the dynamics
results in a huge final-state energy. (Note the different scales). The
naive post-optimization truncation shown here is definitely not a
viable approach to control theory with superpositions of waveforms.
In the bottom panel (c) the control pulse is determined
using the PINV-based approach described in section \ref{sec:method},
applying it to the L-BFGS algorithm \cite{Nocedal06}. The same
subspace as in panel \ref{fig:T0}(b) is used, however, the final energy
$E(\tau)=4.432\cdot 10^{-6}$ is also very close to the exact result.
The initial guess for the control is $\vec{u}_i=0$, and the
convergence parameter is set to $\epsilon=10^{-6}$ for all runs shown in figure
\ref{fig:T0}.

\subsection{Control of an open system at finite temperature}\label{sec:open}

\begin{figure}
\begin{center}
\includegraphics[width=9cm]{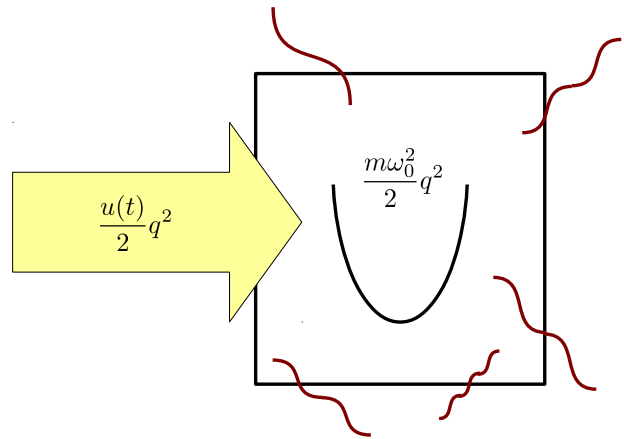}
\caption{(color online) Schematic representation of the system
 under study in section ~\ref{sec:open} and \ref{sec:perform}. A harmonic oscillator is coupled to a thermal
environment, see text for details. The arrow illustrates
 the external parametric control.}\label{fig:system}
\end{center} 
\end{figure}
 
In this section we look again at a harmonic oscillator ($\omega_0=1,
m=1$), but now in contact with an environment in thermal equilibrium
at temperature $T$ as illustrated in figure \ref{fig:system}.

A general description for open classical systems is the Langevin equation \cite{Risken}, which reads for the system under study (with parametric control)
\begin{equation}
\eqalign{\dot{q}(t)=&\frac{p}{m} \cr
\dot{p}(t)=&-m\omega_0^2 q(t)-\int\limits_0^t ds \gamma(t-s)\dot{q}(s) -u(t)q(t) +\xi(t),}
\label{eq:langevin}
\end{equation}
where the stochastic force $\xi(t)$ satisfies the
Fluctuation-Dissipation Theorem, $\langle\xi(t)\xi(t')\rangle=m k_{\rm B}
T\gamma(t-t')$ for $t>t'$. We assume a Drude damping for the
environment associated with a friction-kernel
\begin{equation}
\gamma(t)=\gamma_0 \omega_c \Theta(t)e^{-\omega_c t},
\end{equation}
with damping constant $\gamma_0$ and cut-off frequency $\omega_c$. The
parameters are set to $\gamma_0=0.1$, $\omega_c=10$ and $k_{\rm B} T=1$.

Based on this description valid for single noise realizations $\xi(t)$, the expectation value $\E[\cdot]$ is estimated by taking the average of many realizations. The presented results are achieved with $M=1000$ realizations.

\begin{figure*}
\includegraphics[width=14cm]{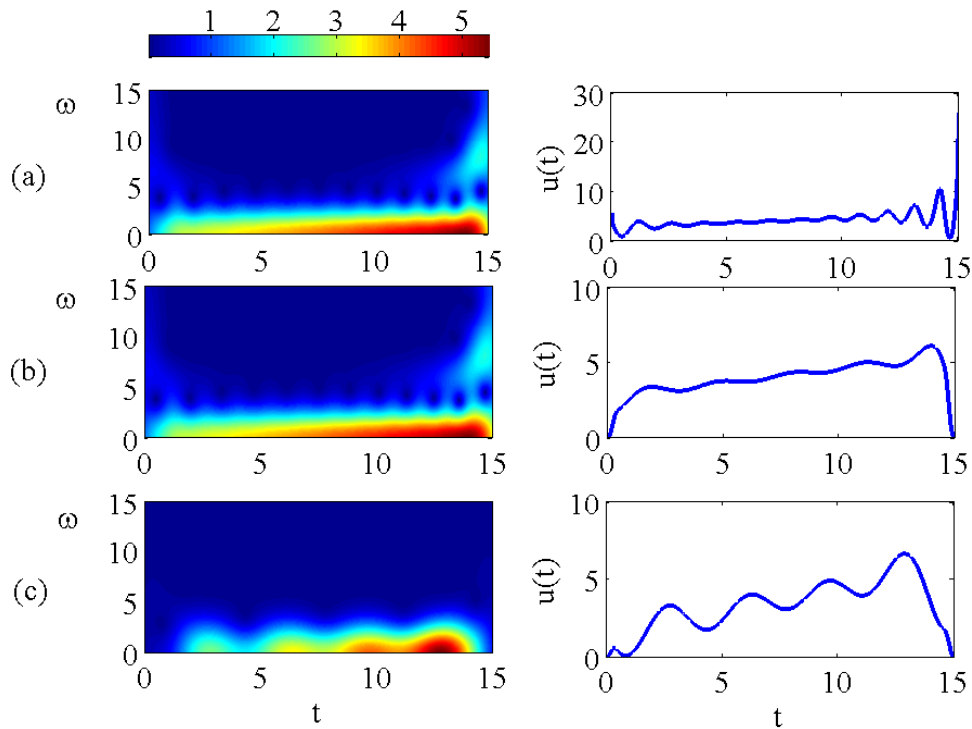}
\caption{(color online) Windowed Fourier Transform of the
control signals (left column) and the time resolved control
 signals (right column) for a parametrically controlled harmonic oscillator in a thermal reservoir at $k_{\rm B}T=1$. 
 Results in (a) correspond to the optimization in full
function space leading to a final-time value  $E(\tau)=0.305$. In (b) this control is
truncated to a $n=12$-dimensional subspace defined in (\ref{eq:basisp}). This control leads to
$E(\tau)=0.500$. In (c) results achieved with the PINV
 method in the same subspace as in (b) are shown producing $E(\tau)=0.323$.
}\label{fig:erg}
\end{figure*}

Initially the system is prepared with Gaussian distributed random
numbers to satisfy $\E[E(0)]=0.5$. In absence of any control the
system would equilibrate to $\E[E]=k_{\rm B}T=1$. The objective here is to
further extract energy from the system towards a low energy state
\begin{equation}
\Phi[q(\tau),p(\tau)]=\E\left[\frac{q^2(\tau)}{2}+\frac{p^2(\tau)}{2}\right].
\end{equation}
First we want to search for a control function in full function space,
as reference solution to which PINV-based results can be compared.  As
the numeric cost is very high to do optimal control in a $N=2250$
dimensional function space, in particular as for every iteration step
(\ref{eq:langevin}) has to be solved for $M=1000$
realizations and averaged, the reference result is computed
more efficiently from a Fokker-Planck equation \cite{Risken}, which
can be mapped to a system of ordinary differential equations in the
case of Gaussian probability densities \cite{Schmidt11,stock14}.

The top panel of figure \ref{fig:erg} shows an optimized full-space
control signal without any constraint obtained using Krotov's algorithm \cite{Krotov96}, as well
as its windowed-Fourier-transformed \cite{Allen77}.  This algorithm is
similar to a gradient search near convergence points, but has the
advantage of better robustness in its dependence on the initial
guess. The resulting control signal is basically a ramp superposed
with periodic oscillations. Shortly before the end time these
oscillations increase markedly in amplitude, leaving the control
signal at a very high value at the final time $\tau=15$. These
features of the unconstrained optimization are probably not
beneficial in an experimental setting. A detailed discussion of this
result is given in \cite{Schmidt12}. With this control the system has
its energy reduced to $E(\tau)=0.305$, well below both its initial and
thermal equilibrium values.

We now turn our attention to optimization constrained to a subspace,
again defined through the modified Fourier basis (\ref{eq:basisp}) with n=12.

The middle panel of figure \ref{fig:erg} shows a naive projection of
the Fokker-Planck result onto this subspace. In the projected control
signal a ramp phase can still be distinguished, but the signal is
markedly altered at the beginning and end of the time
interval. Notably, the high value at the final time is
suppressed. Applying this simple
post-optimization projection leads to a poor result, $E(\tau)=0.500$.

The optimization result obtained with the PINV method, see panel (c) in figure \ref{fig:erg} is
$E(\tau)=0.323$, a result only slightly inferior to that
obtained with the complex signal shown in (a). The shape of the control pulse found with the PINV method is clearly different to the Fokker-Planck result;  the ramp is now superposed with highly developed oscillations. Different control signals yielding comparable values for the
optimization objective are not unusual when a simple system is
controlled for an extended period of time.

When using a gradient method, there is no guarantee that the minimum it
finds is a global one. A good initial guess for the control function
can improve the end result significantly. Starting near a minimum also
reduces the required number of optimization
iterations. Here we use the post-truncated
  Fokker-Planck result (b) as initial guess. The result shown in (c) is obtained by subsequently applying the PINV version of a steepest descent algorithm with
convergence tolerance $\epsilon=10^{-5}$.

As in the deterministic case, we see that the PINV method has
significant advantages over a simple projection of a control signal
obtained in the unconstrained function space. Although the number of
independent variables is drastically reduced in the PINV approach,
optimization results of a quality similar to the unconstrained case
are obtained.

\subsection{Performance considerations
\label{sec:perform}}

\begin{figure*}
\includegraphics[width=12cm]{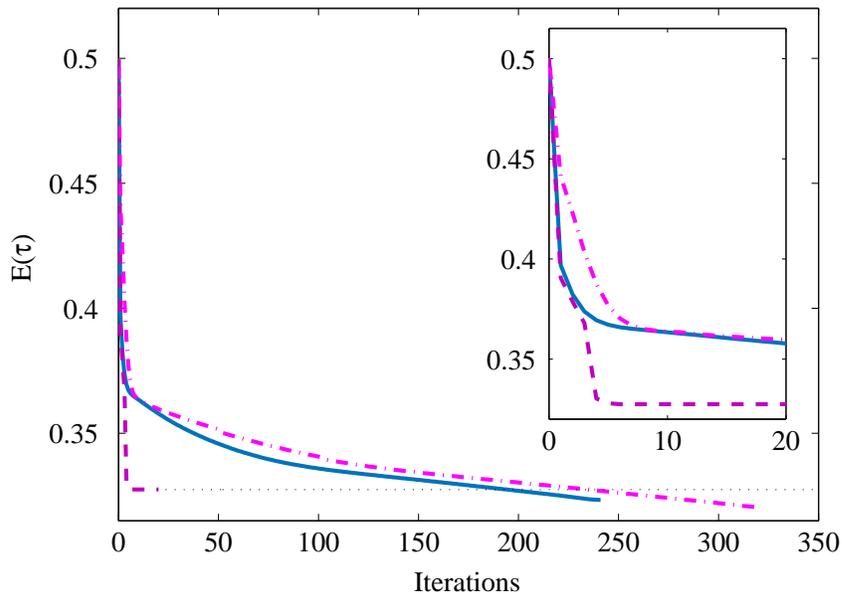}
\caption{(color online) Comparison of convergence behaviour for PINV-projected time-domain
gradient search (solid blue line) and gradient search in coefficient
space (dashed purple line). The latter ends prematurely in a plateau
(see section \ref{sec:perform} for details). This problem is absent if
coefficients of an orthonormal basis are considered (dash-dotted
magenta line). In this case, the iteration progresses in a manner
similar to the time-domain case.}\label{fig:SVD}
\end{figure*}

All results presented so far were obtained using our projection
algorithm in the time domain. An alternative to this approach
\cite{Skinner10} is optimization based on gradients of the objective
with respect to the expansion coefficients $c_l$, which involves the
matrix $B$ as a Jacobian, $\nabla_c J = B\, \nabla_u J$.
This approach has an equivalent representation in the time domain. It
effectively uses the same algorithmic sequence outlined in section
\ref{sec:method}, with one notable change: The projector $B^+B$ is
replaced by the matrix $B^tB$, which differs from the projector unless
an orthonormal basis is chosen.

Figure \ref{fig:SVD} shows a comparison of the two approaches for the
dissipative oscillator discussed in the previous subsection. The
solid line (blue) indicates iterations leading to the result shown in figure \ref{fig:erg}(c). The dashed line (purple)
indicates optimization with the alternative approach of directly using
the coefficient vector. After significant gains for a few iterations,
a plateau is reached; the last iterations shown change the objective
only by about $10^{-6}$. After 20 iterations, standard line
searches fail (the dotted line is continued as a guide to the eye).

This failure seems remarkable since the control landscapes are the
same for both approaches, up to a linear transform mediated by the
matrix $B^t$. We believe the discrepancy between the two approaches is
likely associated with the observation that the linear map described
by $B^tB$ is highly anisotropic; the non-zero singular values of $B$
differ by many orders of magnitude. The observed behaviour is probably
due to rapid convergence along the ``easy axes'' defined by the
anisotropy, followed by extremely slow convergence along the ``hard
axes''.

This interpretation is supported by the following observation: Using
the reduced singular value decomposition \cite{Trefethen97} of the
basis matrix $B=USV^t$, it is possible to find an orthonormal basis of
the subspace through the column vectors of $V$. With $V^t$ substituted
for the matrix $B$, there is no anisotropy; $V V^t$ is identical to
the projector as $B^+B$. An optimization using coefficients for the
orthogonal basis (dash-dotted, magenta) now shows behaviour similar to
the time-domain optimization.

\section{Conclusion and discussion\label{sec:conclusion}}

We have developed a method to adapt standard, gradient-based
techniques to the problem of experimentally realizable control
functions which are characterized by a function space of finite,
typically small dimension. This reduced space is defined as the span
of an arbitrary set of functions. Despite this generality, the
projection of arbitrary functions to the reduced space as well as
transformation between a vector of function values over a discretized
time axis and a coordinate vector in the reduced space can be easily
accomplished by using the Moore-Penrose pseudoinverse.

It is imperative that the iterative computation of optimal control
solutions be performed entirely in the control subspace of
interest. Imposing a single projection on a control solution computed
without constraints yields unacceptable results. It is therefore of
great importance that the transformations and projections implicit in
each iteration step are efficiently computed; this is the case for the
Moore-Penrose pseudoinverse.

The reduced dimensionality of the function space to be searched
typically reduces the numerical cost of optimization, in particular
when using quasi-Newton methods. In addition, we note arbitrary
parametrizations of admissible control fields may introduce
artificial anisotropies to the control problem, potentially leading to
a serious slowdown of gradient methods.

When using projections in the time domain, on the other hand, fairly
rapid convergence can be seen. Similar observations can be made in the
case of parameters which are coordinates of an orthonormal
basis. Iterations already close to an extremum can be accelerated by
using quasi-Newton methods like BFGS, a particular advantage when
high-precision solutions are required, as, e.g., in the context of
quantum information processing.

Our method is easily generalized to non-linear parameterizations of
admissible control functions. In this case, the sets of admissible
control functions are manifolds rather than linear spaces. However,
the present approach can be adapted by computing the projected
gradient in the tangent space at the current iteration point and using
the natural mapping from the parameter space to the manifold.

The class of algorithms studied here may find applications beyond
optimal control theory. Dynamic programming can be thought of as a
formal analogue of optimal control theory, with time replaced by an
arbitrary parameter of a recursion relation formally analogous to an
equation of motion. This opens up the prospect of interdisciplinary
applications or applications outside physics, e.g.,  problems in
economics.

\ack We gratefully acknowledge financial support from the Carl Zeiss Foundation, the Landesgraduiertenstiftung Baden-W\"urttemberg and the DFG through AN336/6-1.

\section*{References}

\bibliographystyle{unsrt}
\bibliography{Bibfile}

\end{document}